\documentclass[draft,12pt]{article}

\usepackage[latin1]{inputenc}
\usepackage{amsmath,amssymb}
\usepackage{latexsym}
\usepackage{amsthm}

\newtheorem{theorem}{Theorem}[section]
\newtheorem{corollary}[theorem]{Corollary}
\newtheorem{lemma}[theorem]{Lemma}
\newtheorem{proposition}[theorem]{Proposition}
\newtheorem{definition}[theorem]{Definition}

\newtheorem{remark}[theorem]{Remark}

\numberwithin{equation}{section}

\def\Var{{\mathop {{\rm Var\, }}}}

\def\square{{\vcenter{\vbox{\hrule height.3pt
        \hbox{\vrule width.3pt height5pt \kern5pt
           \vrule width.3pt}
        \hrule height.3pt}}}}

  \def\sF {{\cal F}}
\def\sG {{\cal G}}

 \def\bQ {{\mathbb Q}}

\def\wt{\widetilde}

\def\ol{\overline}
\def\E{{\mathbb E}}
\def\P{{\mathbb P}}
\def\norm#1{{\Vert #1 \Vert}}

\def\lam{{\lambda}}
\def\angel#1{{\langle #1 \rangle}}
\def\bee{\begin{equation}}
\def\bet{\begin{theorem}}
\def\bep{\begin{proposition}}
\def\bef{\begin{proof}}
\def\bel{\begin{lemma}}
\def\bec{\begin{corollary}}
\def\bed{\begin{definition}}
\def\ber{\begin{remark}}
\def\eee{\end{equation}}
\def\eet{\end{theorem}}
\def\eep{\end{proposition}}
\def\eef{\end{proof}}
\def\eel{\end{lemma}}
\def\eec{\end{corollary}}
\def\eed{\end{definition}}
\def\eer{\end{remark}}

\def\R{{\mathbb R}}
\def\E{{{\mathbb E}\,}}
\def\P{{\mathbb P}}
\def\Q{{\mathbb Q}}

\def\lam{{\lambda}}

\def\al{{\alpha}}

\def\eps{\varepsilon}

\def\angel#1{{\langle#1\rangle}}
\def\norm#1{\Vert #1 \Vert}

 \def\qq {\qquad}

\def\wt{\widetilde}
\def\ol{\overline}

\def\ni{\noindent }
\def\ms{\medskip}

\def\Var{{\mathop {{\rm Var\, }}}}

\def\square{{\vcenter{\vbox{\hrule height.3pt
        \hbox{\vrule width.3pt height5pt \kern5pt
           \vrule width.3pt}
        \hrule height.3pt}}}}

\def\tfrac#1#2{{\textstyle {\frac{#1}{#2}}}}

\def\tlint{{- \kern-0.85em \int \kern-0.2em}}  
\def\dlint{{- \kern-1.05em \int \kern-0.4em}}  

  \def\sF {{\cal F}}
\def\sG {{\cal G}}

 \def\bQ {{\mathbb Q}}

\def\nn{{\nonumber}}

\begin{document}

\title{A stochastic differential equation with a sticky point}
\author{Richard F. Bass}

\date{\today}

\maketitle

\begin{abstract}  
\noindent {\it Abstract:} 
We consider a degenerate stochastic differential equation that has a sticky
point in the Markov process sense. We prove that weak existence and weak
uniqueness hold, but that pathwise uniqueness does not hold nor
does a strong solution exist.

\end{abstract}

\section{Introduction}\label{S:intro}

The one-dimensional stochastic differential equation
\bee\label{intro-E0}
dX_t=\sigma(X_t)\, dW_t
\eee
has been the subject of intensive study for well over half a century.
What can one say about pathwise uniqueness when $\sigma$ is allowed
to be zero at certain points?
Of course,  a large amount is known, but there are
many unanswered questions remaining.

Consider the case where $\sigma(x)=|x|^\al$ for $\al\in (0,1)$.
When $\al \ge 1/2$, it is known there is pathwise uniqueness by the
Yamada-Watanabe criterion (see, e.g., \cite[Theorem 24.4]{stoch}) while if
$\al<1/2$, it is known there are at least two solutions, the zero solution
and one that can be constructed by a non-trivial time change of Brownian
motion.
However, that is not the end of the story. In \cite{xtoal}, it was 
shown 
 that there is in fact pathwise uniqueness when
$\al<1/2$ provided one restricts attention to  the class of solutions that spend zero time at 0. 

This can be better understood by using ideas from Markov process theory.
The continuous strong Markov processes on the real line that are
on natural scale can be characterized by their speed measure. For the
example in the preceding paragraph,  the speed measure $m$ is given by
$$m(dy)=1_{(y\ne 0)} |y|^{-2\al}\, dy+\gamma\delta_0(dy),$$ 
where $\gamma\in [0,\infty]$ and $\delta_0$ is point mass at 0.
When $\gamma=\infty$, we get the 0 solution, or more precisely, the
solution that stays at 0 once it hits 0. If we set $\gamma=0$, we get
the situation considered in \cite{xtoal} where the amount of time
spent at 0 has Lebesgue measure zero, and pathwise uniqueness
holds among such processes.

In this paper we study an even simpler equation:
\bee\label{intro-E1}
dX_t=1_{(X_t\ne 0)}\, dW_t,\qq X_0=0,
\eee
where $W$ is a one-dimensional Brownian motion. 
One solution is $X_t=W_t$, since Brownian motion spends zero time
at 0. Another is the identically 0 solution. 

We take $\gamma\in (0,\infty)$ and consider the class of solutions
to \eqref{intro-E1} which spend a positive amount of time at 0,
with the amount of time parameterized by $\gamma$.
We give a precise description of what we mean by this in Section \ref{S:SMM}.

Representing diffusions on the line as the solutions to  stochastic
differential equations has a long history, going back to It\^o in the 1940's, 
and this paper is a small step in that program. For this reason we characterize
our solutions in terms of occupation times determined by a speed measure. Other
formulations that are purely in terms of stochastic calculus are possible; see
the system \eqref{EPeq1}--\eqref{EPeq2} below.

 We start by proving  weak existence of solutions to \eqref{intro-E1} for
each $\gamma\in (0,\infty)$. We in fact consider a much more 
general situation. We let $m$ be  any measure  that
gives finite positive mass to each open interval and define the
notion of continuous local martingales with speed measure $m$.

 We prove weak uniqueness, or equivalently, uniqueness in law, among
continuous local martingales with speed measure $m$. 
The fact that we have uniqueness in law not only within the class of
strong Markov processes but also within the class of continuous 
local martingales with a given speed measure may be of independent interest.

We then restrict our attention to \eqref{intro-E1} and look at the
class of continuous martingales that solve \eqref{intro-E1} and
at the same time have speed measure $m$, where now
\bee\label{intro-E3}
m(dy)=1_{(y\ne 0)}\, dy+\gamma\delta_0(dy)
\eee
with $\gamma\in (0,\infty)$.

 Even when we fix $\gamma$ and restrict attention to solutions to \eqref{intro-E1}
that have speed measure $m$ given by \eqref{intro-E3},   pathwise uniqueness does not hold.
The proof of this fact is the main result of this paper.
The reader familiar with excursions  will recognize some ideas
from that theory in the proof.

 Finally, we prove that for each $\gamma\in (0,\infty)$, no strong solution
to \eqref{intro-E1} among the class of continuous martingales with
speed measure $m$ given by \eqref{intro-E3} exists. Thus, given $W$, one cannot 
find a continuous martingale  $X$ with speed measure $m$  satisfying
\eqref{intro-E1} such that $X$ is adapted to the
filtration of $W$. A consequence of this is that certain natural
approximations to the solution of \eqref{intro-E1}
do not converge in probability,
although they do converge weakly.

Besides increasing the versatility of \eqref{intro-E0}, one can easily imagine
a practical application of
sticky points. Suppose a corporation has a takeover offer at \$10.
The stock price is then likely to spend a great deal of time 
precisely at \$10
but is not constrained to stay at \$10. Thus \$10 would
be a sticky point for the solution of the stochastic differential equation that describes the stock
price.

Regular continuous strong Markov processes on the line which are 
on natural scale and have speed measure given by \eqref{intro-E3} are
known as sticky Brownian motions. These  were first studied by Feller in the
1950's and It\^o and McKean in the 1960's.

A posthumously published paper by Chitashvili (\cite{Chitashvili})
in 1997, based on a technical report produced in 1988, considered 
processes on the non-negative real line that satisfied the stochastic 
differential equation
\bee\label{one-sided}
dX_t=1_{(X_t\ne 0)}\, dW_t+\theta 1_{(X_t=0)}\, dt, \qq X_t\ge 0, 
\quad X_0=x_0,
\eee
with $\theta\in (0,\infty)$. Chitashvii proved weak uniqueness for the 
pair $(X,W)$ and showed that no strong solution exists.

Warren (see \cite{Warren1} and also \cite{Warren2}) further investigated
solutions to \eqref{one-sided}. The process $X$ is not adapted 
to the filtration generated by $W$ and has some ``extra randomness,''
which Warren characterized.

While this paper was under review, we learned of a preprint by Engelbert and Peskir 
\cite{Engelbert-Peskir} on the subject of sticky
Brownian motions. They considered the system of equations
\begin{align}
dX_t&=1_{(X_t\ne 0)}\, dW_t, \label{EPeq1}\\
1_{(X_t=0)}\, dt&=\frac{1}{\mu}\, d\ell^0_t(X),\label{EPeq2}
\end{align}
where $\mu\in (0,\infty)$ and $\ell^0_t$ is the local time in the
semimartingale sense at 0 of $X$. (Local times in the Markov process sense
can be different in general.) Engelbert and Peskir proved weak uniqueness
of the joint law of $(X,W)$ and proved that no strong solution exists. 
They also considered a one-sided version of this equation, where $X\ge 0$,
and showed that it is equivalent to \eqref{one-sided}. Their results thus
provide a new proof of those of Chitashvili.

It is interesting to compare the system \eqref{EPeq1}--\eqref{EPeq2}
investigated by \cite{Engelbert-Peskir} with the SDE considered in this paper.
Both include the equation \eqref{EPeq1}. In this paper, however, in place of
\eqref{EPeq2} we use a side condition whose origins come from Markov process 
theory, namely:
\begin{align}
X &\mbox{\rm  is a continuous martingale with speed measure }\label{RBeq2}\\
&~~~~~ m(dx)=
dx+\gamma \delta_0(dx),\nn
\end{align}
where $\delta_0$ is point mass at 0 and ``continuous martingale with speed 
measure $m$'' is defined in \eqref{SMM-E1}.  One can show that a solution to the system studied
by \cite{Engelbert-Peskir} is a solution to the formulation considered in this paper and vice versa, and we sketch the argument in 
Remark \ref{comparison}. However, we did not see a way of proving this without 
first proving the uniqueness results of this paper and  using the uniqueness 
results  of \cite{Engelbert-Peskir}.

Other papers that show no strong solution exists for stochastic differential equations that are closely related include 
\cite{Barlow-skew}, \cite{Barlow-LMS}, and \cite{Karatzasetal}.

After a short section of preliminaries, Section \ref{S:prelim}, we
define speed measures for local martingales in Section \ref{S:SMM} and
consider the existence of such local martingales. Section \ref{S:WU}
proves weak uniqueness, while in Section \ref{S:SDE} we prove that
continuous martingales with speed measure $m$ given by \eqref{intro-E3}
satisfy \eqref{intro-E1}.
Sections \ref{S:approx}, \ref{S:est}, and \ref{S:PU} prove that
pathwise uniqueness and strong existence fail. The first of these
sections considers some approximations to a solution to \eqref{intro-E1},
the second proves some needed estimates, and the proof is completed
in the third.

\ni{\bf Acknowledgment.} We would like to thank Prof.\ H.\ Farnsworth
for suggesting a mathematical finance interpretation of a sticky point.

\section{Preliminaries}\label{S:prelim}

For information on martingales and stochastic calculus,
see \cite{stoch}, \cite{KaratzasShreve} or \cite{RevuzYor}.
For background on continuous Markov processes on the line,
see the above references and also
\cite{Ptpde}, \cite{ItoMcKean}, or \cite{Knight}.

We start with an easy lemma concerning continuous local martingales.

\bel\label{prelim-L1} Suppose $X$ is a continuous local martingale
which exits a finite non-empty interval $I$ a.s. If the endpoints
of the interval are $a$ and $b$, $a<x<b$,  and $X_0=x$ a.s.,
then 
$$\E \angel{X}_{\tau_I}=(x-a)(b-x),$$
where $\tau_I$ is the first exit time of $I$ and $\angel{X}_t$ is
the quadratic variation process of $X$.
\eel

\bef Any such local martingale is a time change of a Brownian motion,
at least up until the time of exiting the interval $I$. The result
follows by performing  a change of variables
in the corresponding result for Brownian motion; see, e.g., \cite[Proposition
3.16]{stoch}.
\eef

Let $I$  be a finite non-empty interval with endpoints $a<b$. Each
of the endpoints
may be in $I$ or in $I^c$.
Define $g_I(x,y)$ by
$$g_I(x,y)=\begin{cases} 2(x-a)(b-y)/(b-a),\phantom{\Big]}& a\le x<y\le b;\\
2(y-a)(b-x)/(b-a),& a\le y\le x\le b. \end{cases}$$
Let $m$ be a measure 
such that $m$ gives finite strictly positive measure to every 
finite open interval. Let
$$G_I(x)=\int_I g_I(x,y)\, m(dy).$$

If $X$ is a real-valued process adapted to a filtration $\{\sF_t\}$ 
satisfying the
usual conditions, we let
\bee\label{prelim-E301}
\tau_{I}=\inf\{t>0: X_t\notin I\}.
\eee
When we want to have 
exit times for more than one process at once, we write $\tau_{I}(X)$,
$\tau_{I}(Y)$, etc.
Define
\bee\label{prelim-E302}
T_x=\inf\{t>0: X_t=x\}.
\eee

A continuous strong Markov process $(X,\P^x)$ on the real line is regular
if $\P^x(T_y<\infty)>0$ for each $x$ and $y$. Thus, starting at $x$, there
is positive probability of hitting $y$ for each $x$ and $y$.
A regular continuous strong Markov process
$X$ is on natural scale if whenever $I$ is a  finite non-empty interval
with endpoints $a<b$, then
$$\P^x(X_{\tau_I}=a)=\frac{b-x}{b-a}, \qq \P^x(X_{\tau_I}=b)=\frac{x-a}{b-a}$$
provided $a<x<b$. A continuous regular strong Markov process on the line
on natural scale has speed measure $m$ if for each finite non-empty interval $I$ we have
$$\E^x \tau_I=G_I(x)$$
whenever $x$ is in the interior of $I$.

It is well known that if $(X,\P^x)$ and $(Y,\bQ^x)$ are continuous
regular strong Markov processes on the line on
natural scale with the same speed measure
$m$, then the law of $X$ under $\P^x$ is equal to the law of $Y$ under $\bQ^x$
for each  $x$. 
In addition, $X$ will be a local martingale under $\P^x$ for each
$x$.

Let $W_t$ be a one-dimensional Brownian motion and let $\{L^x_t\}$
be the jointly continuous local times. If we define
\bee\label{cE21}
\al_t=\int L_t^y \, m(dy),
\eee
then $\al_t$ will be continuous and strictly increasing. If we let
$\beta_t$ be the inverse of $\al_t$ and set 
\bee\label{prelim-E1}
X^M_t=x_0+W_{\beta_t},
\eee
then $X^M$ will be a continuous regular strong Markov process on natural
scale with speed measure $m$ starting at $x_0$.   
See the references listed above for a proof, e.g., \cite[Theorem 41.9]{stoch}.
We denote the law of $X^M$ started at $x_0$ by $\P^{x_0}_M$.

If $(\Omega, \sF, \P)$ is a probability space and $\sG$ a $\sigma$-field
contained in $\sF$, a regular
conditional probability $\bQ$ for $\P(\cdot\mid \sG)$ is a map
from $\Omega\times \sF$ to $[0,1]$  such that\\
(1) for each $A\in \sF$, $\bQ(\cdot, A)$ is measurable with respect
to $\sF$;\\
(2) for each $\omega\in \Omega$, $\bQ(\omega,\cdot)$ is a probability
measure on $\sF$;\\
(3) for each $A\in \sF$, $\P(A\mid \sG)(\omega)=\Q(\omega,A)$
for almost every $\omega$.

Regular conditional probabilities do not always exist, but will if
$\Omega$ has sufficient structure; see \cite[Appendix C]{stoch}.
 
The filtration $\{\sF_t\}$ generated by a process $Z$ is the smallest
filtration to which $Z$ is adapted and 
which satisfies the usual conditions.

We use the letter $c$ with or without subscripts to denote finite
positive constants whose value may change from place to place.

\section{Speed measures for local martingales}\label{S:SMM}

Let $a:\R\to \R$ and $b:\R\to \R$ be Borel measurable functions
with $a(x)\le b(x)$ for all $x$. If $S$ is a finite stopping time,
let 
$$\tau^S_{[a,b]}=\inf\{t>S: X_t\notin [a(X_S),b(X_S)]\}.$$

We say a continuous local martingale $X$ started at $x_0$ has speed measure $m$ if
$X_0=x_0$ and 
\bee\label{SMM-E1}
\E[\tau^S_{[a,b]}-S\mid \sF_S]=G_{[a(X_S),b(X_S)]}(X_S), \qq \mbox{\rm a.s.}
\eee
whenever $S$ is a finite stopping time and $a$ and $b$ are as above.

\begin{remark}\label{R-speed}{\rm
 We remark that if $X$ were a strong Markov 
process, then the left hand side of \eqref{SMM-E1} would be equal to 
$\E^{X_S} \tau^0_{[a,b]}$, where $\tau^0_{[a,b]}=\inf\{t\ge 0: X_t\notin [a,b]\}$. Thus the above definition of speed measure for a martingale is a generalization of the one for one-dimensional diffusions on natural scale.
}
\end{remark}

\bet\label{SMM-T1} Let $m$ be a measure that is finite and positive on
every finite open interval.  There exists a continuous local martingale
$X$ with $m$ as its speed measure.
\eet

\begin{proof} Set $X$ equal to $X^M$ as defined in \eqref{prelim-E1}.
We only need show that \eqref{SMM-E1} holds. 
Since $X$ is a Markov
process and has associated with it probabilities $\P^x$ and
shift operators $\theta_t$,  then 
$$\tau^S_{[a,b]}-S=\sigma_{[a(X_0),b(X_0)]}\circ \theta_S,$$ 
where $\sigma_{[a(X_0),b(X_0)]}=\inf\{t>0: X_t\notin [a(X_0),b(X_0)]\}$.
By the strong Markov property, 
\bee\label{SMM-E2}
\E[\tau^S_{[a,b]}-S\mid \sF_S]
=\E^{X_S} \sigma_{[a(X_0), b(X_0)]} \qq \mbox{\rm a.s.}
\eee
For each $y$, $\sigma_{[a(X_0),b(X_0)]}=\tau_{[a(y),b(y)]} $
under $\P^y$, and therefore
$$\E^y \sigma_{[a(X_0),b(X_0)]}=G_{[a(y),b(y)]}(y).$$
Replacing $y$ by $X_S(\omega)$ and substituting in \eqref{SMM-E2}
yields \eqref{SMM-E1}.
\end{proof}

\bet\label{SMM-P1} Let $X$ be any continuous local martingale that has speed
measure $m$ and let $f$ be a non-negative Borel measurable function.
Suppose $X_0=x_0$, a.s.
Let $I=[a,b]$ be a finite interval with $a<b$ such that $m$ does not give positive
mass to  either end point.
Then
\bee\label{SMM-E3m}
\E\int_0^{\tau_I} f(X_s)\, ds=\int_I g_I(x,y) f(y)\, m(dy).
\eee
\eet

\begin{proof} It suffices to suppose that $f$ is continuous and 
equal to 0 at the
boundaries of $I$ and then
to approximate an arbitrary non-negative Borel measurable function by 
continuous functions that are 0 on the boundaries of $I$. 
The main step is to prove
\bee\label{SMM-E4}
\E\int_0^{\tau_I(X)} f(X_s)\, ds=
\E\int_0^{\tau_I(X^M)} f(X^M_s)\, ds.
\eee
Let $\eps>0$. Choose $\delta$ such that
$|f(x)-f(y)|<\eps$ if $|x-y|<\delta$ with $x,y\in I$. 

Set $S_0=0$ and $$S_{i+1}=\inf\{t>S_i: |X_t-X_{S_i}|\ge \delta\}.$$
Then
$$\E\int_0^{\tau_I} f(X_s)\, ds=\E\sum_{i=0}^\infty \int_{S_i\land \tau_I}^{S_{i+1}\land \tau_I} f(X_s)\, ds$$ 
differs by at most $\eps \E\tau_I$ from
\begin{align}
\E\sum_{i=0}^\infty f(X_{S_i\land \tau_I})& (S_{i+1}\land \tau_I-S_i\land \tau_I)\label{SMM-E2a}\\
&=\E\Big[\sum_{i=0}^\infty f(X_{S_i\land \tau_I})\E[S_{i+1}\land \tau_I
-S_i\land \tau_I\mid \sF_{S_i\land \tau_I}]\,\Big].\nn
\end{align}
Let $a(x)=a\lor (x-\delta)$ and $b(x)=b\land (x+\delta)$.
Since $X$ is a continuous local martingale with speed measure $m$,
the last line  in \eqref{SMM-E2a} is equal to
\bee\label{SMM-E2b}
\E\sum_{i=0}^\infty f(X_{S_i\land \tau_I}) 
G_{[a(X_{S_i\land \tau_I}),b( X_{S_i\land \tau_I})]}(X_{S_i\land \tau_I}).
\eee

Because $\E \tau_{[-N,N]}<\infty$ for all $N$, then $X$ is a time change
of a Brownian motion. It follows that  the distribution of $\{X_{S_i\land \tau_I(X)}, i\ge 0\}$ is that of a simple random walk on the lattice $\{x+k\delta\}$
stopped the first time it exits $I$, and thus
is the same as the distribution of $\{X^M_{S_i\land \tau_I(X^M)}, i\ge 0\}$. Therefore
the expression is \eqref{SMM-E2b} is equal to the corresponding expression
with $X$ replaced by $X^M$. This in turns differs by at most
$\E \eps \tau_I(X^M)$ from
$$\E\int_0^{\tau_I(X^M)} f(X^M_s)\, ds.$$
Since $\eps$ is arbitrary, we have
\eqref{SMM-E4}.
Finally, the right hand side of \eqref{SMM-E4}
 is equal to the right hand side of \eqref{SMM-E3m}
by \cite[Corollary IV.2.4]{Ptpde}.
\end{proof}

\section{Uniqueness in law}\label{S:WU}

In this section we show that if $X$ is a continuous local martingale under
$\P$ with speed measure $m$, then $X$ has the same law as $X^M$.
Note that we do not suppose \emph{a priori} that $X$ is a strong
Markov process. We remark that the results of  \cite{ES3} do not apply, since
in that paper a generalization of the system \eqref{EPeq1}--\eqref{EPeq2}
is studied rather than the formulation given by \eqref{EPeq1} together with 
\eqref{RBeq2}.

\bet\label{WU-T1} Suppose $\P$ is a probability measure and $X$
is a continuous local martingale with respect to $\P$. Suppose that
$X$ has speed measure $m$
and $X_0=x_0$ a.s. Then the law of $X$ under $\P$ is equal
to the law of $X^M$ under $\P^{x_0}_M$.
\eet

\bef Let $R>0$ be such that $m(\{-R\})=m(\{R\})=0$ and set $I=[-R,R]$. 
Let $\ol X_t=X_{t\land \tau_I(X)}$ and $\ol X^M_t=X^M_{t\land \tau_I(X^M)}$, 
the processes $X$ and $X^M$ stopped on exiting $I$.
For $f$ bounded and measurable let
$$H_\lam f=\E \int_0^{\tau_I({\ol X})} e^{-\lam t} f({\ol X}_t)\, dt$$
and
$$H_\lam^M f(x)=\E^{x} \int_0^{\tau_I({\ol X}^M)} e^{-\lam t} f({\ol X}^M_t)\, dt$$
for $\lam\ge 0$. Since ${\ol X}$ and ${\ol X}^M$ are stopped at times $\tau_I({\ol X})$
and $\tau_I({\ol X}^M)$, resp.,  we can replace $\tau_I({\ol X})$
and $\tau_I({\ol X}^M)$ by $\infty$ in both of the above integrals 
without affecting $H_\lam$ or $H_\lam^M$ as
long as $f$ is 0 on the boundary of $I$.

Suppose $f(-R)=f(R)=0$. Then $H_\lam^M f(-R)$ and $H_\lam^M f(R)$ are
also 0, since we are working with the stopped process.

We want to show 
\bee\label{WU-E2}
H_\lam f=H_\lam^M f(x_0), \qq  \lam\ge 0.
\eee 
By Theorem \ref{SMM-P1} we know \eqref{WU-E2} holds for $\lam=0$. 
Let $K=\E \tau_I(X)$. We have $\E^{x_0} \tau_I(X^M)=K$ as well
since both $X$ and $X^M$ have speed measure $m$.

Let $\lam=0$ and $\mu\le 1/2K$. Let $t>0$ and let $Y_s={\ol X}_{s+t}$.
Let $\bQ_t$ be a regular conditional probability for $\P(Y\in \cdot\mid
\sF_t)$. 
It is easy to see that for almost every $\omega$, 
$Y$ is a continuous local martingale under $\bQ_t(\omega, \cdot)$ started at ${\ol X}_t$ and $Y$ has speed measure
$m$. Cf.\ \cite[Section I.5]{Ptpde} or \cite{xtoal}. Therefore by Theorem \ref{SMM-P1}
$$\E_{\bQ_t}\int_0^\infty f(Y_s)\, ds=H^M_0 f({\ol X}_t).$$
This can be rewritten as
\bee\label{WU-E31}
\E\Big[\int_0^\infty f({\ol X}_{s+t})\, ds\mid \sF_t\Big]=H^M_0f({\ol X}_t), \qq
\mbox{\rm a.s.}
\eee
as long as $f$ is 0 on the endpoints of $I$.

Therefore, recalling that $\lam=0$, 
\begin{align}
H_\mu H_\lam^M f&=\E\int_0^\infty e^{-\mu t} H_\lam^M f({\ol X}_t)\, dt\label{WU-E3}\\
&=\E\int_0^\infty e^{-\mu t} \E\Big[\int_0^\infty e^{-\lam s}f({\ol X}_{s+t})\, ds\mid
\sF_t\Big]\, dt\nn\\
&=\E\int_0^\infty e^{-\mu t}e^{\lam t}\int_t^\infty e^{-\lam s}f({\ol X}_s)\, ds\, dt\nn\\
&=\E\int_0^\infty \int_0^s e^{-(\mu-\lam)t}\, dt\, e^{-\lam s}f({\ol X}_s)\, ds\nn\\
&=\E\int_0^\infty \frac{1-e^{-(\mu-\lam)s}}{\mu-\lam}e^{-\lam s}f({\ol X}_s)\, ds\nn\\
&=\frac{1}{\mu-\lam}\E\int_0^\infty e^{-\lam s}f({\ol X}_s)\, ds
-\frac{1}{\mu-\lam}\E\int_0^\infty e^{-\mu s} f({\ol X}_s)\, ds.\nn\\
&=\frac{1}{\mu-\lam}H_\lam^M f(x_0)
-\frac{1}{\mu-\lam}\E\int_0^\infty e^{-\mu s} f({\ol X}_s)\, ds.\nn
\end{align} 
We used \eqref{WU-E31} in the second equality.
Rearranging,
\bee\label{SMM-E32}
H_\mu f=H_\lam^M f(x_0)+(\lam-\mu)H_\mu(H_\lam^Mf).
\eee

Since ${\ol X}$ and ${\ol X}^M$ are stopped upon exiting $I$, then $H^M_\lam f=0$
at the endpoints of $I$. We now take \eqref{SMM-E32} with $f$ 
replaced by $H_\lam^M f$, use this  to evaluate the last term in 
\eqref{SMM-E32},  and obtain
$$H_\mu f=H_\lam^Mf({ x}_0)+(\lam-\mu)H^M_\lam(H^M_\lam f)(x_0)+(\lam-\mu)^2 H_\mu(H^M_\lam(H^M_\lam f)).$$
We continue. Since $$|H_\mu g|\le \norm{g} \E\tau_I(X)=\norm{g} K$$
and $$\norm{H^M_\lam g} \le \norm{g} \E \tau_I(X^M)=
\norm{g} K$$
for each bounded $g$, where $\norm{g}$ is the supremum norm of $g$,
we can iterate and get convergence as long as $\mu\le 1/2K$ and obtain
$$H_\mu f=H^M_\lam f(x_0)+\sum_{i=1}^\infty ((\lam-\mu) H^M_\lam)^i H_\lam^M f(x_0).$$

The above also holds when ${\ol X}$ is replaced by ${\ol X}^M$, so that
$$H_\mu^M f(x_0)=H^M_\lam f(x_0)+\sum_{i=1}^\infty ((\lam-\mu) H^M_\lam)^iH_\lam^M  f(x_0).$$
We conclude $H_\mu f=H^M_\mu f(x_0)$ as long as $\mu\le 1/2K$ and
$f$ is 0 on the endpoints of $I$.

This holds for every starting point. If $Y_s={\ol X}_{s+t}$ and $\bQ_t$ is 
a regular conditional probability for the law of $Y_s$ under $\P^x$
given $\sF_t$,
then we asserted above that 
$Y$ is a continuous local martingale started at ${\ol X}_t$ with speed measure $m$
under $\bQ_t(\omega, \cdot)$ 
for almost every $\omega$. We replace $x_0$ by ${\ol X}_t(\omega)$ in the preceding
paragraph and derive
$$\E\Big[\int_0^\infty e^{-\mu s}f({\ol X}_{s+t})\, ds\mid \sF_t\Big]=H^M_\mu f({\ol X}_t),
\qq \mbox{\rm a.s.}$$
if $\mu\le 1/2K$ and $f$ is 0 on the endpoints of $I$.

We now take $\lam=1/2K$ and $\mu\in(1/2K, 2/2K]$. The same argument
as above shows that $H_\mu f=H^M_\mu f(x_0)$ as  long as $f$ is 0 on the 
endpoints of $I$. This is true for every starting point. We continue,
letting $\lam=n/2K$ and using induction, and obtain
$$H_\mu f=H_\mu^M f(x_0)$$
for every $\mu\ge 0$.

Now suppose $f$ is continuous with
 compact support and $R$ is large enough so that $(-R,R)$
contains the support of $f$. We have
that
$$\E\int_0^{\tau_{[-R,R]}({\ol X})} e^{-\mu t} f({\ol X}_t)\, dt
=\E^{x_0}\int_0^{\tau_{[-R,R]}({\ol X}^M)} e^{-\mu t} f({\ol X}^M_t)\, dt$$
for all $\mu> 0$.
This can be rewritten as
\bee\label{WU-E501}
\E\int_0^\infty e^{-\mu t} f({ X}_{t\land \tau_{[-R,R]}(X)})\, dt
=\E^{x_0}\int_0^\infty  e^{-\mu t} f({X}^M_{t\land\tau_{[R,R]}(X^M)})\, dt.
\eee
 If we hold $\mu$ fixed and let $R\to \infty$
in \eqref{WU-E501}, we obtain 
$$\E\int_0^\infty e^{-\mu t} f(X_t)\, dt
=\E^{x_0}\int_0^\infty e^{-\mu t} f(X^M_t)\, dt$$
for all $\mu>0$.
By the uniqueness of the Laplace transform and the continuity of $f, X,$
and $X^M$, 
$$\E f(X_t)=\E^{x_0} f(X^M_t)$$
for all $t$. By a limit argument, this holds whenever $f$ is a bounded Borel
measurable function.

The starting point $x_0$ was arbitrary. Using regular conditional
probabilities as above,
$$\E[f(X_{t+s})\mid \sF_t]=\E^{x_0} [f(X_{t+s}^M)\mid \sF_t].$$
By the Markov property, the right hand side is equal to 
$$\E^{X^M_t} f(X_s)=P_s f(X^M_t),$$
where $P_s$ is the transition probability kernel for $X^M$.

To prove that the finite dimensional distributions of
$X$ and $X^M$ agree, we use induction.
We have
\begin{align*}
\E\prod_{j=1}^{n+1} f_j(X_{t_j})
&=\E_i \prod_{j=1}^{n} f_j(X_{t_j})\E_i[f_{n+1}(X_{t_{n+1}})\mid \sF_{t_n}]\\
&=\E_i \prod_{j=1}^{n} f_j(X_{t_j})P_{t_{n+1}-t_n}f_{n+1}(X_{t_n}).
\end{align*}
We use the induction hypothesis to see that this is equal to 
$$\E^{x_0} \prod_{j=1}^{n} f_j(X^M_{t_j})P_{t_{n+1}-t_n}f_{n+1}(X^M_{t_n}).$$
We then use the Markov property to see that this in turn is equal to 
$$\E^{x_0}\prod_{j=1}^{n+1} f_j(X^M_{t_j}).$$

Since $X$ and $X^M$ have continuous paths and the same finite dimensional
distributions, they have the same law.
\eef

\section{The stochastic differential equation}\label{S:SDE}

We now discuss the particular stochastic differential equation we want
our martingales to solve. We specialize to the following
speed measure. Let $\gamma\in (0,\infty)$ and
let
\bee\label{SDE-E31}
m(dx)= dx+\gamma\delta_0(dx),
\eee
where $\delta_0$ is point mass at 0.

We consider the stochastic differential equation 
\bee\label{SMM-E200}
X_t=x_0+\int_0^t 1_{(X_s\ne 0)}\, dW_s.
\eee
A triple $(X,W,\P)$ is a  weak solution to \eqref{SMM-E200} 
with $X$ starting at $x_0$ if $\P$ is a probability measure, 
there exists a filtration $\{\sF_t\}$ satisfying the usual conditions,
 $W$ is a Brownian motion under $\P$
with respect to $\{\sF_t\}$, 
and $X$ is a continuous martingale adapted to $\{\sF_t\}$ with $X_0=x_0$ and satisfying  \eqref{SMM-E200}.

We now show that any martingale with $X_0=x_0$ a.s.\ that has speed measure $m$ is the first
element of a triple that is a weak solution to \eqref{SMM-E200}.
Although $X$ has the same law as $X^M$ started at $x_0$, here we only have
one probability measure and we cannot assert that $X$ is a strong Markov 
process.
We point out that  \cite[Theorem 5.18]{ES3} does not apply here, since they study a generalization of the system \eqref{EPeq1}--\eqref{EPeq2}, and we do not know at this stage that this formulation is equivalent to the one used here.

\bet\label{SMM-T3} Let $\P$ be a probability measure on a space that supports a Brownian motion and let
$X$ be a continuous martingale which has speed measure $m$ with $X_0=x_0$ a.s. 
Then there
exists a Brownian motion $W$ such that $(X,W,\P)$ is a weak
solution to \eqref{SMM-E200} with $X$ starting at $x_0$.
Moreover
\bee\label{sde-E323}
X_t=x_0+\int_0^t 1_{(X_s\ne 0)}\, dX_s.
\eee
\eet

\begin{proof}
Let
$$W'_t=\int_0^t  1_{(X_s\ne 0)}\, dX_s.$$
Hence
$$d\angel{W'}_t=1_{(X_t\ne 0)}\, d\angel{X}_t.$$

Let $0<\eta<\delta$. Let $S_0=\inf\{t: |X_t|\ge \delta\}$,
$T_i=\inf\{t>S_i: |X_t|\le \eta\}$, and $S_{i+1}=\inf\{t>T_i:
|X_t|\ge \delta\}$ for $i=0,1,\ldots$.

The speed measure of
$X$ is equal to $m$, which in turn is equal 
to  Lebesgue measure on $\R\setminus\{0\}$,  
hence $X$ has the same law as $X^M$ by Theorem \ref{WU-T1}. Since $X^M$ behaves 
like a Brownian motion when it is away from zero, 
we conclude $1_{[S_i,T_i]}\,d\angel{X}_t=1_{[S_i,T_i]}\, dt$. 

Thus
for each $N$,
$$\int_0^t  1_{\cup_{i=0}^N [S_i,T_i]}(s)\, d\angel{X}_s
=\int_0^t  1_{\cup_{i=0}^N [S_i,T_i]}(s)\, ds.$$
Letting $N\to \infty$, then $\eta\to 0$, and finally $\delta\to \infty$,
we obtain
$$\int_0^t 1_{(X_s\ne 0)}\, d\angel{X}_s=\int_0^t 1_{(X_s\ne 0)}\, ds.$$

Let $V_t$ be an independent Brownian motion and let
$$W''_t=\int_0^t 1_{(X_s=0)}\, dV_s.$$ 
Let $W_t=W'_t+W''_t$. Clearly $W'$ and $W''$ are orthogonal martingales, so
$$d\angel{W}_t=d\angel{W'}_t+d\angel{W''}_t
=1_{(X_t\ne 0)}\, dt+1_{(X_t=0)}\, dt=dt.$$
By L\'evy's theorem (see \cite[Theorem 12.1]{stoch}), $W$ is a Brownian motion.

If 
$$M_t=\int_0^t 1_{(X_s=0)}\, dX_s,$$
by the occupation times formula (\cite[Corollary VI.1.6]{RevuzYor}),
$$\angel{M}_t=\int_0^t 1_{(X_s=0)}\, d\angel{X}_s
=\int 1_{\{0\}}(x) \ell^x_t(X)\, dx=0$$
for all $t$,
where $\{\ell^x_t(X)\}$ are the local times of $X$ in the semimartingale sense.
This implies that $M_t$
is identically zero, and hence $X_t=W'_t$.

Using the definition of $W$, we deduce
\bee\label{sde-E5.95}
1_{(X_t\ne 0)} \, dW_t=1_{(X_t\ne 0)}\, dX_t=dW'_t=dX_t,
\eee
as required.
\end{proof}

We now show weak uniqueness, that is, if $(X,W,\P)$ and $(\wt X,\wt W,\wt \P)$
are two weak solutions to \eqref{SMM-E200} with $X$ and $\wt X$ starting 
at $x_0$  and in addition $X$ and $\wt X$ have
speed measure $m$, then the joint law of
$(X,W)$ under $\P$ equals the joint law of $(\wt X,\wt W)$ under $\wt \P$.
This holds even though $W$ will not in general be adapted
to the filtration of $X$.
We know that the law of $X$ under $\P$ equals the law of $\wt X$ under
$\wt \P$ and also that
the law of $W$ under $\P$ equals the law of $\wt W$ under
$\wt \P$, but the issue here is the joint law.
Cf.\ \cite{Cherny}. See also \cite{Engelbert-Peskir}.

\bet\label{WU-T21} Suppose $(X,W,\P)$ and $(\wt X,\wt W,\wt \P)$
are two weak solutions to \eqref{SMM-E200}  with $X_0=\wt X_0=x_0$ 
and that $X$ and $\wt X$
are both continuous martingales with speed measure $m$.
 Then the joint law of
$(X,W)$ under $\P$ equals the joint law of $(\wt X,\wt W)$ under $\wt \P$.
\eet

\bef 
Recall the construction of $X^M$ from Section \ref{S:prelim}.
With $U_t$ a Brownian motion  with jointly continuous local times $\{L^x_t\}$
and $m$ given by \eqref{SDE-E31}, we define $\al_t$ by \eqref{cE21}, let
$\beta_t$ be the right continuous inverse of $\al_t$, and let $X^M_t=x_0+U_{\beta_t}$. 
Since $m$ is greater than or equal to  Lebesgue measure  but is finite on every
finite interval, we see that $\al_t$ is strictly increasing, continuous,
 and $\lim_{t \to \infty} \al_t=\infty$. It follows that $\beta_t$ is 
continuous and tends
to infinity almost surely as $t\to \infty$. 

Given any stochastic process $\{N_t, t\ge 0\}$, let $\sF^N_\infty$ be the 
$\sigma$-field generated by the collection of random variables $\{N_t, t\ge 0\}$ together with  the null sets.

We have $\beta_t=\angel{X^M}_t$ and $U_t=X^M_{\al_t}-x_0$. 
Since $\beta_t$ is measurable with respect to $\sF^{X^M}_\infty$
for each $t$, then $\al_t$ is also, and hence so is $U_t$. In fact, we 
can give a recipe to construct
a Borel measurable map $F:C[0,\infty)\to C[0,\infty)$ such that $U=F(X^M)$.
 Note 
also that $X^M_t$ is measurable with respect
to $\sF^U_\infty$ for each $t$ and there exists a Borel measurable map
$G:C[0,\infty)\to C[0,\infty)$ such that $X^M=G(U)$. In addition observe that $\angel{X^M}_\infty=\infty$ a.s.

Since $X$ and $X^M$ have the
same law, then $\angel{X}_\infty=\infty$ a.s. If $Z_t$ is a Brownian 
motion with $X_t=x_0+Z(\zeta_t)$ for a continuous increasing process $\zeta$, then
$\zeta_t=\angel{X}_t$ is measurable with respect to $\sF^X_\infty$,
 its inverse $\rho_t$ is also, and therefore $Z_t=X_{\rho_t}-x_0$ is as 
well. 
Moreover the recipe for constructing $Z$ from $X$ is exactly the same as the one for
constructing $U$ from $X^M$, that is, 
 $Z=F(X)$.
Since $X$ and $X^M$ have the same law, then  the joint law of  $(X,Z)$ is equal to the joint law of $(X^M,U)$.
We can therefore conclude that $X$ is measurable with respect to 
$\sF^Z_\infty$ and $X=G(Z)$.

Let $$Y_t=\int_0^t 1_{(X_s=0)}\, dW_s.$$
Then $Y$ is a martingale with $$\angel{Y}_t=\int_0^t 1_{(X_s=0)}\, ds
=t-\angel{X}_t.$$
Observe that $\angel{X,Y}_t=\int_0^t 1_{(X_s\ne 0)}1_{(X_s=0)}\, ds=0$.
By a theorem of Knight (see \cite{Knight2} or \cite{RevuzYor}), there exists
a two-dimensional process $V=(V_1,V_2)$ such that $V$ is a two-dimensional
Brownian motion under $\P$ and 
$$(X_t,Y_t)=(x_0+V_1(\angel{X}_t),V_2(\angel{Y}_t), \qq \mbox{\rm a.s.}$$
(It turns out that $\angel{Y}_\infty=\infty$, but that is not needed in Knight's
theorem.)

By the third paragraph of this proof,  $X_t=x_0+V_1(\angel{X}_t)$ implies that $X_t$ is measurable with
respect to $\sF^{V_1}_\infty$, and in fact $X=G(V_1)$.
  Since $\angel{Y}_t=t-\angel{X}_t$,
then 
$(X_t,Y_t)$ is measurable with respect to $\sF^V_\infty$ for each $t$ and there
exists a Borel measurable map $H: C([0,\infty), \R^2)\to C([0,\infty), \R^2)$,
where $C([0,\infty), \R^2)$ is  the space of continuous functions from
$[0,\infty)$ to $\R^2$, and $(X,Y)=H(V)$. 
Thus $(X,Y)$ is the image under $H$ of a two-dimensional Brownian motion. If 
$(\wt X, \wt W, \wt \P)$ is another weak solution, then we can define $\wt Y$ 
analogously and find a two-dimensional Brownian motion $\wt V$ such that
$(\wt X, \wt Y)=H(\wt V)$. The key point is that the same $H$ can be used.
We conclude that the law of $(X,Y)$ is 
uniquely determined. 
Since $$(X,W)=(X,X+Y-x_0),$$ this proves that the joint law
of $(X,W)$ is uniquely determined.
\eef

\begin{remark}\label{comparison}{\rm
In Section \ref{S:prelim} we constructed the continuous strong Markov process $(X^M,\P_M^x)$
and we now know that $X$ started at $x_0$ is equal in law to $X^M$ under
$\P^{x_0}_M$. We pointed out  in Remark \ref{R-speed} that in the strong Markov case the notion of speed
measure for a martingale reduces to that of speed measure for a one dimensional diffusion. In \cite{Engelbert-Peskir} it is shown that the solution to the
system \eqref{EPeq1}--\eqref{EPeq2} is unique in law and thus the solution
started at $x_0$ is equal in law to that of a diffusion on $\R$ started at $x_0$; let $\wt m$ be
the speed measure for this strong Markov process. Thus to show the equivalence
of the system \eqref{EPeq1}--\eqref{EPeq2} to the one given by \eqref{EPeq1} and
\eqref{RBeq2},
it suffices to show that $\wt m=m$
if and only if \eqref{EPeq2} holds, where $m$ is
given by \eqref{SDE-E31} and $\gamma=1/\mu$. Clearly both $\wt m$ and $m$ are equal to Lebesgue measure on $\R\setminus 
\{0\}$, so it suffices to compare the atoms of $\wt m$ and $m$ at 0.

Suppose \eqref{EPeq2} holds and $\gamma=1/\mu$. Let $A_t=\int_0^t 1_{\{0\}}(X_s)\, ds$. Thus 
\eqref{EPeq2} asserts that $A_t=\frac{1}{\mu}\ell_t^0$. Let $I=[a,b]=[-1,1]$,
$x_0=0$, and $\tau_I$ the first time that $X$ leaves the interval 
$I$. Setting $t=\tau_I$ and taking expectations starting from 0, we 
have
$$\E^0 A_{\tau_I}=\frac1{\mu} \E^0 \ell^0_{\tau_I}.$$
Since $\ell^0_t$ is the increasing part of the submartingale $|X_t-x_0|-|x_0|$
and $X_{\tau_I}$ is equal to either 1 or $-1$, the right hand side is equal to 
$$\frac{1}{\mu} \E^0|X_{\tau_I}|=\frac1{\mu}.$$
On the other hand, by \cite[(IV.2.11)]{Ptpde}, 
$$\E^0 A_{\tau_I}=\int_{-1}^1 g_I(0,y) 1_{\{0\}}(y)\, \wt m(dy)
=\wt m(\{0\}).$$
Thus $\wt m=m$ if $\gamma=1/\mu$.

Now suppose we have a solution to the pair \eqref{EPeq1} and \eqref{RBeq2}
and $\gamma=1/\mu$; we will
show \eqref{EPeq2} holds. Let $R>0$, $I=[-R,R]$, and $\tau_I$ the first exit
time from $I$. Set $B_t=\frac{1}{\mu} \ell^0_t$.  For any $x\in I$, we have 
by \cite[(IV.2.11)]{Ptpde} that
\bee\label{Rc1}
\E^x A_{\tau_I}=\int_{-1}^1 g_I(x,y)1_{\{0\}}(y)\, m(dy)
=\gamma g_I(x,0).
\eee 
Taking expectations,
\bee\label{Rc2}
 \E^x B_{\tau_I}=\frac{1}{\mu}\E^x[\,|X_{\tau_I}-x|-|x|\,].
\eee
Since $X$ is a time change of a Brownian motion that exits $I$ a.s., the distribution of $X_{\tau_I}$ started at $x$ is the same as that of a Brownian motion
started at $x$ upon exiting $I$. A simple computation shows that the
right hand side of \eqref{Rc2} agrees with the right hand side of \eqref{Rc1}.
By the strong Markov property,
$$\E^0[A_{\tau_I}-A_{\tau_I\land t}\mid \sF_t]=
\E^{X_t} A_{\tau_I}=\E^{X_t}B_{\tau_I}
=\E^0[B_{\tau_I}-B_{\tau_I\land t}\mid \sF_t]$$
almost surely on the set $(t\le \tau_I)$.
Observe that 
if $U_t=\E^0[A_{\tau_I}-A_{\tau_I\land t}\mid \sF_t]$, then we can write
$$U_t
=\E^0[A_{\tau_I}-A_{\tau_I\land t}\mid \sF_t]
=\E^0[A_{\tau_I}\mid \sF_t] -A_{\tau_I\land t}$$
and $$U_t=\E^0[B_{\tau_I}-B_{\tau_I\land t}\mid \sF_t]
=\E^0[B_{\tau_I}\mid \sF_t] -B_{\tau_I\land t}$$
for $t\le \tau_I$.
This expresses the supermartingale  $U$ as  a martingale minus an increasing process
in two different ways.
By the uniqueness of the Doob decomposition for supermartingales, we conclude
$A_{\tau_I\land t}=B_{\tau_I\land t}$
for $t\le \tau_I$. Since $R$ is arbitrary, this establishes \eqref{EPeq2}.
(The argument that the potential of an increasing process determines the process
is well known.)
}\end{remark}

\begin{remark}\label{eitherpaper}{\rm
In the remainder of the paper we prove that there does not exist a 
strong solution to the pair \eqref{EPeq1} and \eqref{RBeq2} nor does pathwise
uniqueness hold. In \cite{Engelbert-Peskir}, the authors prove that there is no
strong solution to the pair \eqref{EPeq1} and \eqref{EPeq2} and that pathwise
uniqueness does not hold. Since we now know there is an equivalence between the pair
\eqref{EPeq1} and \eqref{RBeq2} and the pair \eqref{EPeq1} and \eqref{EPeq2},
one could at this point use the argument of \cite{Engelbert-Peskir} in place of  the argument of this paper. 
Alternatively, in the paper of \cite{Engelbert-Peskir} one could use our argument in place of theirs to establish the non-existence of a strong solution and 
that pathwise uniqueness does not hold.
}\end{remark}

\section{Approximating processes}\label{S:approx}

Let $\wt W$ be a Brownian motion adapted to a filtration $\{\sF_t, t\ge 0\}$, let $\eps\le \gamma$,  and let $X_t^\eps$ be the solution to 
\bee\label{approx-E671}
dX^\eps_t=\sigma_\eps(X_t^\eps)\, d\wt W_t, \qq X^\eps_0=x_0,
\eee
where
$$\sigma_\eps(x)=\begin{cases} 1,& |x|>\eps;\\
\sqrt{\eps/\gamma},& |x|\le \eps.\end{cases}$$
 For
each $x_0$ the solution to the stochastic differential
equation is pathwise unique by \cite{LeGall} or \cite{Nakao}. 
We also know that if $\P^x_\eps$ is the law of $X^\eps$ starting
from $x$, then $(X^\eps, \P^x_\eps)$ is a continuous regular
strong Markov process on natural scale. The speed measure
 of $X^\eps$ will be 
$$m_\eps(dy)=dy+\frac{\gamma}{\eps} 1_{[-\eps,\eps]}(y)\, dy.$$
Let $Y^\eps$ be the solution to 
\bee\label{approx-E672}
dY^{\eps}_t=\sigma_{2\eps}(Y^\eps_t)\, d\wt W_t, \qq Y^\eps_0=x_0.
\eee

Since $\sigma_\eps\le 1$, then $d\angel{X^\eps}_t\le dt$. By
the Burkholder-Davis-Gundy inequalities 
(see, e.g., \cite[Section 12.5]{stoch}), 
\bee\label{c6.2A}
\E|X^\eps_t-X^\eps_s|^{2p}\le c|t-s|^p
\eee
for each $p\ge 1$, where the constant $c$ depends on $p$. It follows
(for example, by Theorems 8.1 and 32.1 of  \cite{stoch}) that the law of $X^\eps$ is tight in $C[0,t_0]$ for each
$t_0$. The same is of course true
for $Y^\eps$ and $\wt W$, and so the triple $(X^\eps,Y^\eps,\wt W)$
is tight in $(C[0,t_0])^3$ for each $t_0>0$.

Let $P_t^\eps$ be the transition probabilities for the Markov process
$X^\eps$.
Let $C_0$ be the set of continuous functions on $\R$ that vanish at
infinity and let 
$$L=\{f\in C_0: |f(x)-f(y)|\le |x-y|, x,y\in \R\},$$
the set of Lipschitz functions with Lipschitz constant 1 that vanish  at
infinity.

One of the main results of \cite{Lsg} (see Theorem 4.2) is that $P_t^\eps$
maps $L$ into $L$ for each $t$ and each $\eps<1$.

\bet\label{approx-T101}
If $f\in C_0$, then $P_t^\eps f$ converges uniformly for each $t\ge 0$.
If we denote the limit by $P_tf$, then $\{P_t\}$ is a family of transition
probabilities for a continuous regular strong Markov process $(X, \P^x)$ on natural
scale with speed measure given by \eqref{SDE-E31}.
For each $x$, $\P^x_{\eps}$ converges weakly to $\P^x$ with respect to 
$C[0,N]$ for each $N$.
\eet

\bef \ni\emph{Step 1.}
Let $\{g_j\}$ be a countable collection of $C^2$ functions in  $L$ with
compact support such that the set of finite linear combinations of
elements of $\{g_j\}$ is dense in $C_0$ with respect to the supremum
norm.

Let $\eps_n$ be a sequence converging to 0.
Suppose  $g_j$ has support contained in $[-K,K]$ with $K>1$. Since 
$X_t^\eps$ is a Brownian motion outside $[-1,1]$, if $|x|>2K$, then
$$|P_t^\eps g_j(x)|=|\E^x g_j(X_t^\eps)|\le \norm{g_j}\,
\P^x(|X^\eps|\mbox{ hits $|x|/2$ before time } t),$$
which tends to 0 uniformly over $\eps<1$ as $|x|\to \infty$.
Here $\norm{g_j}$ is the supremum norm of $g_j$.
By the equicontinuity of the $P_t^{\eps} g_j$, using the diagonalization
method there exists a subsequence, which we continue to denote
by $\eps_n$,  such that $P_t^{\eps_n} g_j$ 
converges uniformly on $\R$ for every rational $t\ge 0$
and every $j$. We denote the limit by $P_tg_j$. 

Since $g_j\in C^2$,
\begin{align*}
P_t ^\eps g_j(x)-P_s^\eps g_j(x)&=\E^x g_j(X_t^\eps)-\E^x g_j(X_s^\eps)\\
&=\E^x \int_s^t \sigma_\eps(X_r^\eps)g'_j(X_r^\eps)\, d\wt W_r+\tfrac12 \E^x \int_s^t \sigma_\eps(X_r^\eps)^2g''_j(X_r^\eps)\, dr\\
&=\tfrac12 \E^x \int_s^t \sigma_\eps(X_r^\eps)^2g''_j(X_r^\eps)\, dr,
\end{align*}
where we used Ito's formula. Since $\sigma_\eps$ is bounded by 1, we obtain
$$|P_t^\eps g_j(x)-P_s^\eps g_j(x)|\le c_j |t-s|,$$
where the constant $c_j$ depends on $g_j$. With this fact, we can deduce
that $P_t^{\eps_n} g_j$ converges uniformly in $C_0$ for every $t\ge 0$.
We again call the limit $P_t g_j$. Since linear combinations of the 
$g_j$'s are dense in $C_0$, we conclude that $P_t^{\eps_n} g$ converges
uniformly to a limit, which we call $P_tg$, whenever $g\in C_0$. We note
that $P_t$ maps $C_0$ into $C_0$.

\ni\emph{Step 2.} Each $X_t^\eps$ is a Markov process, so $P_s^\eps(P_t^\eps g)=P_{s+t}^\eps g$.
By the uniform convergence and equicontinuity and the fact
that $P_s^\eps$ is a contraction, we see that $P_s(P_tg)=P_{s+t}g$
whenever $g\in C_0$.

Let $s_1<s_2<\cdots s_j$ and let $f_1, \ldots f_j$ be elements of $L$.
Define inductively $g_j=f_j$, $g_{j-1}=f_{j-1}(P_{s_j-s_{j-1}}g_j)$, 
$g_{j-2}=f_{j-2}(P_{s_{j-1}-s_{j-2}}g_{j-1})$, and so on. Define $g_j^\eps$
analogously where we replace $P_t$ by $P_t^\eps$. By the Markov property
applied repeatedly,
$$\E^x[f_1(X^\eps_{s_1})\cdots f_j(X^\eps_{s_j})] =P^\eps_{s_1} g_1^\eps(x).$$

Suppose $x$ is fixed for the moment
and let $f_1, \cdots, f_j\in L$. Suppose there is a subsequence $\eps_{n'}$ of $\eps_n$ such
that $X^{\eps_{n'}}$ converges weakly, say to $X$,
 and let $\P'$ be the limit law with corresponding expectation $\E'$. Using the
uniform convergence, the equicontinuity, and the fact that $P_t^\eps$ maps $L$ into
$L$, we obtain
\bee\label{Approx-E31}
\E'[f_1(X_{s_1})\cdots f_j(X_{s_j})] =P_{s_1} g_1(x).
\eee

We can conclude several things from this. First, since the limit is the same no
matter what subsequence $\{\eps_{n'}\}$ we use, then the full sequence $\P^x_{\eps_n}$
converges weakly.
This holds for each starting point $x$. 

Secondly, if we denote the weak limit of the $\P^x_{\eps_n}$ by $\P^x$,
then \eqref{Approx-E31} holds with $\E'$ replaced by $\E^x$.
From this
we deduce that $(X,\P^x)$  is a Markov process with transition semigroup
given by $P_t$. 

Thirdly, since $\P^x$ is the weak limit of probabilities on $C[0,\infty)$, 
we conclude that $X$ under $\P^x$ has continuous paths for each $x$.

\ni \emph{Step 3.} Since $P_t$ maps $C_0$ into $C_0$ and $P_tf(x)=\E^x f(X_t)\to f(x)$
by the continuity of paths if $f\in C_0$, we conclude by \cite[Theorem 20.9]{stoch} that
$(X,\P^x)$ is in fact a strong Markov process.

Suppose $f_1, \ldots, f_j$ are in $L$ and $s_1<s_2<\cdots < s_j<t<u$.
Since $X_t^\eps$ is a martingale,
$$\E^x_\eps\Big[X_u^\eps\prod_{i=1}^j f_i(X_{s_i}^\eps) \Big]=
\E^x\Big[X_t^\eps\prod_{i=1}^j f_i(X_{s_i}^\eps) \Big].$$
Moreover, $X_t^\eps$ and $X_u^\eps$ are uniformly integrable due to 
\eqref{c6.2A}.
Passing to the limit along the sequence $\eps_n$, we have the equality
with $X^\eps$ replaced by $X$ and $\E^x_\eps$ replaced by $\E^x$. 
Since the collection of  random variables of the form $\prod_i f_i(X_{s_i})$ generate
$\sigma(X_r; r\le t)$, it follows that $X$ is a martingale under
$\P^x$ for each $x$.

\ni\emph{Step 4.} 
 Let $\delta, \eta>0$. Let $I=[q,r]$ and 
$I^*=[q-\delta,r+\delta]$.
In this step we show that
\bee\label{approx-E102}
\E \tau_I(X)= \int_I g_I(0,y)\, m(dy).
\eee

First we obtain a uniform bound on $\tau_{I^*}(X^\eps)$.
If
$A^\eps_t=t\land \tau_{I^*}(X^\eps)$, then 
$$\E[A^\eps_\infty-A^\eps_t\mid \sF_t]=\E^{X_t^\eps} A^\eps_\infty
\le \sup_x \E^x \tau_{I^*}(X^\eps).$$
The last term is equal to 
$$\sup_x \int_{I^*} g_{I^*}(x,y) \Big(1+\frac{\gamma}{\eps}1_{I^*}(y)\Big)\, dy.$$
A simple calculation shows that this is bounded by $$c(r-q+2\delta)^2+c\gamma (r-q+2\delta),$$
where $c$ does not depend on $r, q, \delta$,  or $\eps$.
By Theorem I.6.10 of \cite{pta}, we then deduce that 
$$ \E \tau_{I^*}(X^\eps)^2
=\E (A^\eps_\infty)^2<c<\infty,$$
where $c$ does not depend on $\eps$.
By Chebyshev's inequality, for each $t$,
$$\P(\tau_{I^*}(X^\eps)\ge t)\le c/t^2.$$

Next we obtain an upper bound on $\E \tau_I(X)$ in terms of $g_{I^*}$.
We have 
\begin{align*}
\P(\tau_I(X)>t)&=\P(\sup_{s\le t} |X_s|\le r,\inf_{s\le t}|X_s|\ge q)\\
&\le 
\limsup_{\eps_n\to 0} \P(\sup_{s\le t}|X_s^{\eps_n}|\le r,\inf_{s\le t}|X_s|^\eps\ge q)\\
&\le \limsup_{\eps_n\to 0} \P(\tau_{I^*}(X^{\eps_n})>t)\le c/t^2.
\end{align*}

Choose $u_0$ such that 
$$\int_{u_0}^\infty \P(\tau_I(X)>t)\, dt<\eta, \qq
\int_{u_0}^\infty \P(\tau_{I^*}(X^{\eps_n})>t)\, dt<\eta$$
for each $\eps_n$.

Let $f$ and $g$ be continuous functions taking values in $[0,1]$ 
such that $f$ is equal
to $1$ on $(-\infty,r]$ and $0$ on $[r+\delta, \infty)$
and $g$ is equal to 1 on $[q,\infty)$ and 0 on $(-\infty,q-\delta]$.
We have
\begin{align*}
\P(\sup_{s\le t}|X_s|\le r,\inf_{s\le t}|X_s\ge q)
&\le \E[ f(\sup_{s\le t}|X_s|)g(\inf_{s\le t}|X_s|)]\\
&=\lim_{\eps_n\to 0} \E[ f(\sup_{s\le t}|X_s^{\eps_n}|)g(\inf_{s\le t}|X_s^{\eps_n}|)].
\end{align*}
Then
{\allowdisplaybreaks
\begin{align*}
\int_0^{u_0} \P(\tau_I(X)>t)\, dt&=\int_0^{u_0} \P(\sup_{s\le t}|X_s|\le r,\inf_{s\le t}|X_s|\ge q)\, dt\\
&\le \int_0^{u_0}\E[ f(\sup_{s\le t}|X_s|)g(\inf_{s\le t}|X_s|)]\, dt\\
&=\int_0^{u_0} \lim_{\eps_n\to 0} \E [f(\sup_{s\le t}|X_s^{\eps_n}|)g(\inf_{s\le t}|X^{\eps_n}_s|)]\, dt\\
&=\lim_{\eps_n\to 0} \int_0^{u_0} \E [f(\sup_{s\le t}|X_s^{\eps_n}|)g(\inf_{s\le t}|X^{\eps_n}_s|)]\, dt\\
&\le \limsup_{\eps_n\to 0}\int_0^{u_0}\P(\sup_{s\le t}|X_s^{\eps_n}|\le r+\delta,\inf_{s\le t}|X_s|\ge q-\delta)\, dt\\
&\le \limsup_{\eps_n\to 0} \int_0^{u_0}\P(\tau_{I^*}(X^{\eps_n})\ge t)\, dt\\
&\le \limsup_{\eps_n\to 0} \E \tau_{I^*}(X^{\eps_n}).
\end{align*}
}

Hence
$$\E \tau_I(X)\le \int_0^{u_0}\P(\tau_I(X)>t)\, dt+\eta
\le \limsup_{\eps_n\to 0} \E \tau_{I^*}(X^{\eps_n})+\eta.$$
We now use the fact that $\eta$ is arbitrary and let $\eta\to 0$.
Then
\begin{align*}
\E \tau_I(X)&\le \limsup_{\eps_n\to 0} \E \tau_{I^*}(X^{\eps_n})\\
&=\limsup_{\eps_n\to 0} \int_{I^*} g_{I^*}(0,y)\Big(1+\frac{\gamma}{\eps}1_{[-\eps,\eps]}(y)\Big)\, dy\\
&=\int_{I^*} g_{I^*}(0,y)m(dy).
\end{align*}

We next use
 the joint continuity of $g_{[-a,a]}(x,y)$ in the variables $a,x$ and $y$.
Letting $\delta\to 0$, we obtain
$$\E \tau_I(X)\le \int_I g_I(0,y)\, m(dy).$$
The lower bound for $\E \tau_I(X)$ is done similarly,
and we obtain
\eqref{approx-E102}.

\ni\emph{Step 5.}
Next we show that $X$ is a regular strong Markov process. This means
that if $x\ne y$, $\P^x(X_t=y\mbox{ for some }t)>0$. To show this,
assume without loss of generality that $y<x$. Suppose $X$ starting from
$x$ does not hit $y$ with positive probability. 
Let $z=x+4|x-y|$. Since $\E^x\tau_{[y,z]}<\infty$, then with probability
one $X$ hits $z$ and does so before hitting $y$. 
Hence $T_z=\tau_{[y,z]}<\infty$ a.s. Choose $t$ large
so that $\P^x(\tau_{[y,z]}>t)<1/16$. 
By the optional stopping theorem,
$$\E^x X_{T_z\land t}\ge z\P^x(T_z\le t)+y\P^x(T_z> t)
=z-(z-y)\P^x(T_z>t).$$
By our choice of $z$, this is greater than $x$, which contradicts that
$X$ is a martingale. Hence $X$ must hit $y$ with positive probability.

Therefore $X$ is a regular continuous strong Markov process on the 
real line. Since it is a martingale, it is on natural scale. Since
its speed measure is the same as that of $X^M$ by \eqref{approx-E102}, we conclude from
\cite[Theorem IV.2.5]{Ptpde} that $X$ and $X^M$ have the same law.
In particular, $X$ is a martingale with speed measure $m$.

\ni\emph{Step 6.} Since we obtain the same limit law no matter what sequence
$\eps_n$ we started with, the full sequence $P_t^\eps$ converges to $P_t$
and $\P^x_\eps$ converges weakly to $\P^x$ for each $x$.

All of the above applies equally well to $Y$ and its transition probabilities
and laws.
\eef

Recall that the sequence $(X^\eps, Y^\eps, \wt W)$ is tight with respect to 
$(C[0,N])^3$ for each $N$.
Take a subsequence $(X^{\eps_n}, Y^{\eps_n}, \wt W)$ that
converges weakly, say to the triple  $(X,Y,W)$, with respect to $(C[0,N])^3$
 for
each $N$. 
The last task of this section is to prove that $X$ and $Y$  satisfy \eqref{SMM-E200}.

\bet\label{approx-T55}
$(X,W)$ and $(Y,W)$ each satisfy \eqref{SMM-E200}.
\eet

\bef
We prove this for $X$ as the proof for $Y$ is exactly the same.
Clearly $W$ is a Brownian motion.
Fix $N$. 
We will first show 
\bee\label{Approx-E64}
\int_0^t 1_{(X_s\ne 0)}\, dX_s=\int_0^t 1_{(X_s\ne 0)}\, dW_s
\eee
 if $t\le N.$ 

Let $\delta>0$ and let $g$ be a continuous function taking values
in $[0,1]$ such that $g(x)=0$ if $|x|<\delta$ and $g(x)=1$ if $|x|\ge
2\delta$. 
Since $g$ is bounded and continuous and $(X^{\eps_n},\wt W)$ converges
weakly to $(X,W)$, then $(X^{\eps_n}, \wt W, g(X^{\eps_n}))$ converges
weakly to $(X,W, g(X))$.  Moreover, since $g$ is 0 on $(-\delta, \delta)$, then
\bee\label{CE1}
\int_0^t g(X^{\eps_n}_s)\, d\wt W_s= \int_0^t g(X^{\eps_n}_s)\, dX^{\eps_n}_s
\eee
for ${\eps_n}$ small enough.

By Theorem 2.2 of \cite{Kurtz-Protter}, we have
$$\Big(\int_0^t g(X^{\eps_n}_s)\, d\wt W_s, \int_0^t g(X^{\eps_n}_s)\, dX^{\eps_n}_s\Big)$$
converges weakly to 
$$\Big(\int_0^t g(X_s)\, dW_s, \int_0^t g(X_s)\, dX_s\Big).$$
Then
\begin{align*}\E&\arctan\Big(\Big|\int_0^t g(X_s)\, dW_s- \int_0^t g(X_s)\, dX_s\Big|\Big)\\
&=\lim_{n\to \infty}
\E\arctan\Big(\Big|\int_0^t g(X^{\eps_n}_s)\, d\wt W_s- \int_0^t g(X^{\eps_n}_s)\, dX^{\eps_n}_s\Big|\Big)=0,
\end{align*}
or
$$\int_0^t g(X_s)\, dW_s= \int_0^t g(X_s)\, dX_s, \qq\mbox{\rm a.s.}$$
Letting $\delta\to 0$ proves \eqref{Approx-E64}.

We know $$X^M_t=\int_0^t 1_{(X^M_s\ne 0)}\, dX^M_s.$$
Since $X^M$ and $X$ have the same law, the same is true if we replace $X^M$ by $X$.
Combining with \eqref{Approx-E64} proves
\eqref{SMM-E200}.
\eef

\section{Some estimates}\label{S:est}

Let $$j^\eps(s)=\begin{cases} 1,&|X^\eps_s|\in [-\eps,\eps] \mbox{ or } 
|Y^\eps_s|\in [-2\eps,2\eps] \mbox{ or both};\\
0,& \mbox{otherwise}.\end{cases}$$ Let 
$$J^\eps_t=\int_0^t j^\eps_s\, ds.$$
Set $$Z_t^\eps=X_t^\eps-Y^\eps_t,$$
suppose $Z_0^\eps=0$,
and define $\psi_\eps(x,y)=\sigma_\eps(x)-\sigma_{2\eps}(y)$.
Then
$$dZ^\eps_t=\psi_\eps(X^\eps_t, Y^\eps_t)\, d\wt W_t.$$

Let 
\begin{align}
S_1&=\inf\{t: |Z_t^\eps|\ge 6\eps\},\label{est-E21}\\
T_i&=\inf\{t\ge S_i: |Z_t^\eps|\notin [4\eps,b]\},\nn\\
S_{i+1}&=\inf\{t\ge T_i: |Z_t^\eps|\ge 6\eps\}, \qq\mbox{\rm  and}\nn\\
U_b&=\inf\{t: |Z_t^\eps|=b\}.\nn
\end{align}

\bep\label{est-P1}
For each $n$,
$$\P(S_n< U_b)\le \Big(1-\frac{2\eps}{b}\Big)^n.$$
\eep

\begin{proof}
Since $X^\eps$ is a recurrent diffusion, $\int_0^t 1_{[-\eps,\eps]}(X_s^\eps)
\, ds$ tends to infinity a.s.\ as   $t\to \infty$.
When $x\in [-\eps,\eps]$, then $|\psi_{\eps}(x,y)|\ge c\eps$, and we
conclude that $\angel{Z^\eps}_t\to \infty$ as $t\to \infty$. 

Let $\{\sF_t\}$ be the filtration generated by $\wt W$.
$Z^\eps_{t+S_n}-Z^\eps_{S_n}$ is a martingale started at 0 with respect
to the regular conditional probability for the law of $(X^\eps_{t+S_n}, Y^\eps_{t+S_n})$ given $\sF_{S_n}$.
The conditional probability that  
 it hits $4\eps$ before $b$ 
if $Z^\eps_{S_n}=6\eps$ is
the same as the conditional probability it hits $-4\eps$ before $-b$ if
$Z_{S_n}^\eps=-6\eps$ and is equal to 
$$\frac{b-6\eps}{b-4\eps}\le 1-\frac{2\eps}{b}.$$
Since this is independent of $\omega$, we have 
$$\P\Big(|Z^\eps_{t+S_n}-Z^\eps_{S_n}|\mbox{ hits }4\eps \mbox{ before hitting } b
\mid \sF_{S_n}\Big)\le 1-\frac{2\eps}{b}.$$

Let $V_n=\inf\{t> S_n: |Z_t^\eps|=b\}$.
Then 
\begin{align*}
\P(S_{n+1}<U_b)&\le \P(S_n<U_b, T_{n+1}<V_n)\\
&=\E[\P(T_{n+1}<V_n\mid \sF_{S_n}); S_n<U_b]\\
&\le \Big(1-\frac{2\eps}{b}\Big)\P(S_n<U_b).
\end{align*}
Our result follows by induction.
\end{proof}

\bep\label{est-P2} There exists a constant $c_1$ such that
$$\E J^\eps_{T_n}\le c_1n\eps$$
for each $n$.
\eep

\begin{proof} For $t$ between times $S_n$ and $T_n$ we know that
$|Z_t^\eps|$ lies between
$4\eps$ and $b$. Then at least one of  $X^\eps_t\notin [-\eps,\eps]$
and  $Y^\eps_t\notin [-2\eps,2\eps]$ holds. 
 If exactly one holds, then  $|\psi_\eps(X_t^\eps,Y_t^\eps)|\ge 
1-\sqrt{2\eps/\gamma}\ge 1/2$ if $\eps$ is small enough. If
both hold, we can only say that $d\angel{Z^\eps}_t\ge 0$. In any case, 
$$d\angel{Z^\eps}_t\ge \tfrac14\,dJ^\eps_t$$
for $S_n\le t\le T_n$.

$Z_t^\eps$ is a martingale, and by Lemma \ref{prelim-L1} 
and an argument  using regular conditional probabilities similar to those
we have done earlier,
\bee\label{est-E301}
\E[J^\eps_{T_n}-J^\eps_{S_n}]\le 4\E[\angel{Z^\eps}_{T_n}
-\angel{Z^\eps}_{S_n}]\le 4(b-6\eps)(2\eps)=c\eps.
\eee

Between times $T_n$ and $S_{n+1}$ it is possible that
$\psi_\eps(X_t^\eps,Y_t^\eps)$ can be 0 or it can be
larger than $c\sqrt{\eps/\gamma}$. However if either $X_t^\eps\in [-\eps,\eps]$
or $Y_t^\eps\in [-2\eps,2\eps]$, then
$\psi_\eps(X_t^\eps,Y_t^\eps)\ge c\sqrt{\eps/\gamma}$.
Thus 
$$d\angel{Z^\eps}_t\ge c\eps\,dJ^\eps_t$$
for $T_n\le t\le S_{n+1}$. By Lemma \ref{prelim-L1}
\bee\label{est-E302}
\E[J^\eps_{S_{n+1}}-J^\eps_{T_n}]\le c\eps^{-1}\E[\angel{Z^\eps}_{S_{n+1}}
-\angel{Z^\eps}_{T_n}]\le c\eps^{-1}(2\eps)(10\eps)=c\eps.
\eee

Summing each of \eqref{est-E301} and \eqref{est-E302} over $j$ from 1 to $n$
and combining  yields the proposition.
\end{proof}

\bep\label{est-P3} Let $K>0$ and $\eta>0$. There exists $R$ depending on 
$K$ and $\eta$ such that
$$\P(J^\eps_{\tau_{[-R,R]}(X^\eps)}<K)\le \eta, \qq \eps\le 1/2.$$
\eep

\begin{proof} Fix $\eps\le 1/2$. We will see that our estimates are
independent of $\eps$. Note
$$J^\eps_t\ge H_t=\int_0^t 1_{[-\eps,\eps]}(X_s^\eps)\,ds.$$
Therefore to prove the proposition it is enough to prove that
$$\P^0_\eps(H_{\tau_{[-R,R]}(X^\eps)}<K)\le \eta$$
if $R$ is large enough.

Let $I=[-1,1]$.
We have
$$\E^0_\eps H_{\tau_I(X^\eps)}\ge \int_{-1}^1 g_{I} (0,y)\frac{\gamma}{\eps} 1_{[-\eps,\eps]}(y)
\, dy\ge c_1.$$
On the other hand, for any $x\in I$,
$$\E^x_0 H_{\tau_{I}(X^\eps)}=\int_I g_I(x,y) \frac{\gamma}{\eps}1_{[-\eps,\eps]}(y) \,
dy\le c_2.$$
Combining this with
$$\E^0_\eps[H_{\tau_I(X^\eps)}-H_t\mid \sF_t]\le \E^{X_t^\eps}_\eps H_{\tau_{I}(X^\eps)}$$
and Theorem I.6.10 of \cite{pta} (with $B=c_2$ there), we see that
$$\E H_{\tau_I(X^\eps)}^2\le c_3.$$

Let $\al_0=0$, $\beta_i=\inf\{t>\al_i: |X_t^\eps|=1\}$ and
$\al_{i+1}=\inf\{t>\beta_i: X_t^\eps=0\}.$
Since $X_t^\eps$ is a recurrent diffusion, each $\al_i$ is finite a.s.\ and
$\beta_i\to \infty$ as $i\to \infty$. 
Let $V_i=H_{\beta_i}-H_{\al_i}$. By the strong Markov property,
under $\P^0_\eps$ the $V_i$ are i.i.d.\ random variables with mean larger
than $c_1$ and variance bounded
by $c_4$, where $c_1$ and $c_4$ do not depend on $\eps$ as long as $\eps<1/2$.
Then
\begin{align*}
\P^0_\eps\Big(\sum_{i=1}^k V_i\le c_1k/2\Big)&\le \P^0_\eps\Big(\sum_{i=1}^k (V_i-\E V_i)\ge c_1k/2\Big)\\
&\le \frac{\Var(\sum_{i=1}^k V_i)}{(c_1k/2)^2}\\
&\le 4c_4/c_1^2k.
\end{align*}
Taking $k$ large enough, we see that
$$\P^0_\eps\Big(\sum_{i=1}^k V_i\le K\Big)\le \eta/2.$$

Using the fact that $X_t^\eps$ is a martingale, starting at 1, the probability
of hitting $R$ before hitting 0 is  $1/R$. Using the strong Markov property,
the probability of $|X|$ having no more than  $k$ downcrossings of $[0,1]$ before
 exiting
$[-R,R]$ is bounded by
$$1-\Big(1-\frac{1}{R}\Big)^k.$$
If we choose $R$ large enough, this last quantity will be less than
$\eta/2$. Thus, except for an event of probability at most $\eta$, $X_t^\eps$
will exit $[-1,1]$ and return to 0 at least $k$ times before
exiting $[-R,R]$ and the total amount of time spent in $[-\eps,\eps]$
before exiting $[-R,R]$ will be at least $K$.
\end{proof}

\bep\label{est-P4} Let $\eta>0, R>0,$ and $I=[-R,R]$. There exists
$t_0$ depending on $R$ and $\eta$  such that
$$\P^0_\eps(\tau_I(X^\eps)>t_0)\le \eta, \qq \eps\le 1/2.$$
\eep

\bef
If $\eps\le 1$,
$$\E^0_\eps \tau_R(X^\eps)=\int_I g_I(x,y)\, m_\eps(dy).$$
A calculation shows this is bounded by $cR^2+cR$,
where
$c$ does not depend on $\eps$ or $R$. Applying Chebyshev's inequality,
$$\P^0_\eps(\tau_I(X^\eps)>t_0)\le \frac{\E^0_\eps\tau_I(X^\eps)}{t_0},$$
which is bounded by $\eta$ if $t_0\ge c(R^2+R)/\eta$.
\eef

\section{Pathwise uniqueness fails}\label{S:PU}

We continue the notation of Section \ref{S:est}.
The strategy of proving that pathwise uniqueness does not hold owes a great
deal to \cite{Barlow-LMS}.

\bet\label{PU-T1} There exist three processes $X,Y$, and $W$ and
a probability measure $\P$ such that $W$ is a Brownian motion
under $\P$, $X$ and $Y$ are continuous martingales under $\P$
with speed measure $m$ starting at 0, \eqref{SMM-E200} holds for $X$, 
\eqref{SMM-E200} holds when $X$ is replaced by $Y$, and
$\P(X_t\ne Y_t\mbox{ for some }t>0)>0$.
\eet

\bef
Let $(X^\eps, Y^\eps, \wt W)$ be defined as in \eqref{approx-E671} and
\eqref{approx-E672} and
choose a sequence $\eps_n$ decreasing to 0 such that 
the triple converges weakly on $C[0,N]\times C[0,N]\times C[0,N]$
for each $N$. By Theorems \ref{approx-T101} and \ref{approx-T55}, the weak limit,
$(X,Y,W)$ is such that $X$ and $Y$ are continuous martingales with
speed measure $m$, $W$ is a Brownian motion, and \eqref{SMM-E200}
holds for $X$ and also when $X$ is replaced by $Y$.

Let $b=1$ and
let $S_n$, $T_n$, and $U_b$ be defined by \eqref{est-E21}.
Let $A_1(\eps,n)$ be the event where $T_n< U_b$.
By Proposition \ref{est-P1}
$$\P(A_1(\eps,n))= \P(S_n< U_b)\le \Big(1-\frac{2\eps}{b}\Big)^n.$$
Choose $n\ge \beta/\eps$, where $\beta$ is large enough so that
the right hand side is less than $1/5$ for all $\eps$ sufficiently
small. 

By Proposition \ref{est-P2},
$$\E J_{T_n}^\eps\le c_1n \eps=c_1\beta.$$
By Chebyshev's inequality,
$$\P(J^\eps_{T_n}\ge 5c_1\beta)\le \P(J_{T_n}^\eps\ge 5\E J_{T_n}^\eps)
\le 1/5.$$
Let $A_2(\eps, n)$ be the event where $J^\eps_{T_n}\ge 5c_1\beta$.

Take $K=10c_1\beta$. By Proposition \ref{est-P3},
there exists $R$ such that
$$\P(J^\eps_{\tau_{[-R,R]}(X^\eps)}<K)\le 1/5.$$
Let $A_3(\eps,R,K)$ be the event where $J^\eps_{\tau_{[R,R]}(X^\eps)}<  K$.

Choose $t_0$ using Proposition \ref{est-P4}, so that except for an
event of probability $1/5$ we have  $\tau_{[-R,R]}(X^\eps)\le t_0$.
Let $A_4(\eps,R, t_0)$ be the event where $\tau_{[-R,R]}(X^\eps)\le t_0$.

Let 
$$B(\eps)=(A_1(\eps,n)\cup A_2(\eps,n)\cup A_3(\eps,R,K)\cup 
A_4(\eps, R, t_0))^c.$$
Note $\P(B(\eps))\ge 1/5$.

Suppose we are on the event $B(\eps)$. We have
$$J^\eps_{T_n}\le 5c_1\beta< K\le J^\eps_{\tau_{[-R,R]}(X^\eps)}.$$
We conclude that  $T_n< \tau_{[-R,R]}(X^\eps)$. 
Therefore, on the event $B(\eps)$,   we see that
$T_n$ has occurred  before time $t_0$. We also know that $U_b$ has occurred before time $t_0$.
Hence, on $B(\eps)$,
$$\P(\sup_{s\le t_0} |Z_s^\eps|\ge b)\ge 1/5.$$

Since $Z^\eps=X^\eps-Y^\eps$ converges weakly to $X-Y$,
then with probability at least $1/5$, we have that $\sup_{s\le t_0} |Z_s|\ge b/2$.
This implies that $X_t\ne Y_t$ for some $t$, or pathwise uniqueness
does not hold.
\eef

We also can conclude that strong existence does not hold. 
The argument we use is similar to ones given in \cite{Cherny}, \cite{Engelbert},
and \cite{Kurtz}.

\bet\label{PU-T2} Let $W$ be a Brownian motion. There does
not exist a continuous martingale $X$ 
starting at 0
with speed measure $m$ 
such that \eqref{SMM-E200} holds and such that $X$ is measurable
with respect to the filtration of $W$.
\eet

\bef
Let $W$ be a Brownian motion and suppose
there did exist such  a process $X$.
Then there is a measurable map
$F:C[0,\infty)\to C[0,\infty)$ such that $X=F(W)$.

Suppose $Y$ is any other continuous martingale with
speed measure $m$  satisfying \eqref{SMM-E200}. Then
by Theorem \ref{WU-T1}, the law of $Y$ equals the law of $X$, and
by Theorem \ref{WU-T21}, the joint law of $(Y,W)$ is equal
to the joint law of $(X,W)$. Therefore $Y$ also satisfies 
$Y=F(W)$, and we get pathwise uniqueness since $X=F(W)=Y$.
However, we know pathwise uniqueness does not hold. We conclude
that no such $X$ can exist, that is,  strong existence does not hold.
\eef


\medskip

\ni {\bf Richard F. Bass}\\
Department of Mathematics\\
University of Connecticut \\
Storrs, CT 06269-3009, USA\\
{\tt r.bass@uconn.edu}
\ms

\end{document}